\documentstyle[amssymb,amstex,graphicx,graphics,12pt,a4wide]{article}
\newtheorem{theorem}{Theorem}

\newtheorem{proposition}{Proposition}

\newtheorem{remark}{Remark}

\newcommand{\Facet}{\operatorname{Facet}}

\newcommand{\Cells}{\operatorname{Cells}}
\newcommand{\TypicalCell}{\operatorname{TCell}}

\newcommand{\Vol}{\operatorname{Vol}}

\newcommand{\ero}{\operatorname{ero}}

\begin{document}

\title{\textbf{Shape-Driven Nested Markov Tessellations}}
\author{Tomasz Schreiber (Toru\'n) and Christoph Th\"ale (Osnabr\"uck)}
\date{}
\maketitle

\begin{abstract}
A new and rather broad class of stationary (i.e. stochastically translation invariant) random tessellations of the $d$-dimensional Euclidean space is introduced, which are called shape-driven nested Markov tessellations. Locally, these tessellations are constructed by means of a spatio-temporal random recursive split dynamics governed by a family of Markovian split kernel, generalizing thereby the -- by now classical -- construction of iteration stable random tessellations. By providing an explicit global construction of the tessellations, it is shown that under suitable assumptions on the split kernels (shape-driven), there exists a unique time-consistent whole-space tessellation-valued Markov process of stationary random tessellations compatible with the given split kernels. Beside the existence and uniqueness result, the typical cell and some aspects of the first-order geometry of these tessellations are in the focus of our discussion.
\end{abstract}
\begin{flushleft}\footnotesize
\end{flushleft}

\begin{center}
\parbox{12cm}{\textit{In December 2010, at the age of $35$, my colleaugue and friend Tomasz died, when our paper was already near its final form. I have tried to complete the text in his spirit.}}
\end{center}

\section{Introduction}

Modern stochastic geometry has a rapidly growing demand for non-trivial and flexible, yet mathematically tractable models for spatial random structures of the $d$-dimensional Euclidean space ${\Bbb R}^d$, $d\geq 2$. The need for such new model arises in particular in the theory of random tessellations, as they are widely used for modeling structures ranging from classical geology to the more recent telecommunication science (cf. \cite{SAPPL, SKM}). Stationary random tessellations that are stable under the operation of iteration -- so-called STIT tessellations -- were introduced recently by Mecke, Nagel and Wei\ss\ \cite{MNW, MNW2, NW05} and they may serve as a new mathematical reference model beside the classical Poisson hyperplane or Poisson-Voronoi tessellations. On the other hand, the first-order geometry of these tessellations is already determined by two parameters, the surface density and the so-called directional distribution. Moreover, STIT tessellations have Poisson-typical cells, which is to say that the distribution of the typical cell of a STIT tessellations coincides with that of a Poisson hyperplane tessellation having the same parameters. Thus, the size (measured for example by some of the intrinsic volumes) of the cells of STIT tessellations is rather inhomogeneous and varies from arbitrarily small to rather large. However, tessellations with such behavior are often not suitable for practical applications. For these reasons, it is desirable to introduce more flexibility into the model. One could think for example of a hard-core condition, which excludes in some sense 'small' cells or of a condition that keeps the size of large cells under control. In the Poisson-Voronoi reps. the Poisson-Delaunay case a related theory was developed in \cite{D1, D2}. A first attempt in the case of iteration stable and related models was recently made by Cowan \cite{Cowan}, where a number of random cell splitting mechanisms of \textit{finite} volumes are discussed, which are mostly \textit{independent} of the cell geometry (so-called geometry independent apportionments). However, neither the problem of existence and uniqueness of a related whole-space tessellation nor main properties of such tessellations (such as the typical cell distribution or scaling relations for example) have been established in that paper and remained open questions.\\ In the present work, we start by describing a general class of what we call nested tessellations. Locally, their construction can be interpreted as outcome of a random process of cell division, which may roughly be described as follows: Each cell is provided with a random life-time, which is related to the geometry of the cell, to its environment and even to the running construction time. Once this life-times is expired, a random hyperplane hitting the cell is chosen according to some distribution, which is allowed to \textit{depend} on the cell's geometry, its environment and also on time. It is introduced in the cell, cut-off by the cell's boundary and divides it into two polyhedral sub-cells. The whole construction is continued independently and recursively in both of the sub-cells, in such a way that the respectively new introduced hyperplanes are also cut-off by the boundaries of their mother-cells. Formally, the construction is described by a family of Markovian split-kernel and uses the general framework of pure jump Markov processes. In \cite{SLDP} it was shown that the same informal description can also be provided for the \textit{infinite volume} construction. Thus, each general nested tessellation admits a unique family of split kernels. Naturally, it is of particular interest if also the converse is also true, i.e. whether a given family of split kernels determines a nested tessellation and, moreover, if this tessellation is unique and 'well behaved'. It seems that problem, in its general form is very difficult to handle.\\ It is the main purpose of the present paper to prove a result of this kind for the class of what we call shape-driven nested Markov tessellations. They are nested tessellations in the above sense, whose characterizing split kernels satisfy a number of regularity conditions. Beside the existence and uniqueness theorem, we also investigate the geometry all these so-constructed tessellations have in common. In the planar case, the most explicit results are available, whereas for $d\geq 3$, mean values depend on an increasing number of parameters, determined by the precise cell splitting mechanism. For this reason and because of the rapidly increasing geometric complexity of the structures, most of our mean value formulas will be restricted to the practically relevant cases $d=2$ and $d=3$. Moreover, we discuss in detail the geometry of the typical cell, which eventually leads to a simulation method via Monte-Carlo techniques.\\ \\ The paper is structured as follows: In Section \ref{secSTIT} we recall in some detail the -- by now well known-- description of STIT tessellations and figure-out the main construction principle underlying the structure of general nested tessellations in terms of their so-called split kernels. After that, in Section \ref{secMarkovTess}, the class of shape-driven Markov tessellations is introduced by assuming certain, yet flexible, regularity conditions of the split kernels and the principal existence and uniqueness result is formulated. Section \ref{secGeometry} is focused on general mean values for these tessellations in terms of parameters determined by the cell spliting mechanism. The argument of our main theorem as well as a global construction of the tessellations under consideration is the content of Section \ref{secProof}.

\section{Iteration Stable and General Nested Tessellations}\label{secSTIT}

 We consider a translation-invariant measure $\Lambda$ on the set ${\cal H}$ of 
 $(d-1)$-dimensional (affine) hyperplanes in ${\Bbb R}^d$ with $d\geq 2$ and note that under the standard
 identification of a hyperplane $H = \{ x \in {\Bbb R}^d,\; \langle x, u \rangle = r \},\;
 u \in {\Bbb S}_{d-1},\; r \geq 0$ with the pair $(u,r) \in {\Bbb S}_{d-1} \times {\Bbb R}_+,$
 the translational invariance leads to the decomposition $\Lambda = \rho \, {\cal R} \otimes \ell_+$, where
 $\ell_+$ is the Lebesgue measure on ${\Bbb R}_+,$ ${\cal R}$ is an even
 probability measure on the unit sphere ${\Bbb S}_{d-1}$ and
 $\rho \in {\Bbb R}_+$ is the {\bf surface density} of $\Lambda.$ In the particular case where
 $\rho = 1$ and where ${\cal R}$ is the uniform distribution on ${\Bbb S}_{d-1}$, the
 corresponding measure is isometry invariant and is denoted by
 $\Lambda_{\rm iso}$ in this paper.  Often, $\Lambda_{\rm iso}$ is referred to as the invariant hyperplane measure.\\ \\ An iteration stable random tessellation $Y_{\Lambda}(t)$ with {\bf driving measure} $\Lambda$ and time parameter $t>0$ is constructed according to the following rules, originally due to Mecke, Nagel and Weiss, see \cite{MNW, MNW2, NW05}, and hence referred
 to as the \textbf{MNW-construction} in the sequel. We begin at the time $0$ with a compact and convex subset $W\subset{\Bbb R}^d$ (called window here) having nonempty interior, regard it as the initial cell
 for the tessellation under construction and assign to it an exponentially distributed
 random lifetime with parameter $\Lambda([W])$, where $$[W] := \{ H \in {\cal H},\; H
 \cap W \neq \emptyset \}$$ stands for the set of all hyperplanes hitting $W$.
 Upon expiry of its lifetime, the cell $W$ dies
 and splits into two sub-cells $W^+$ and $W^-$ separated by a hyperplane in $[W]$, chosen
 according to the law \begin{equation}\Lambda(\cdot)/\Lambda([W]).\label{SplitMeasureSTIT}\end{equation} The resulting new cells $W^+$ and $W^-$ are again assigned independent exponentially distributed random lifetimes with respective parameters $\Lambda([W^+])$ and $\Lambda([W^-])$ (whence geometrically smaller cells live stochastically longer)
 and the entire construction continues recursively, until the deterministic time threshold
 $t>0$ is reached. The cell-separating $(d-1)$-dimensional facets arising in subsequent
 splits are usually referred to as \textbf{maximal polytopes} (which are in this paper always $(d-1)$-dimensional) or I-segments in the planar case $d=2$, as assuming shapes similar to the literal {\it I}. We will also speak about maximal polygons in the spatial case $d=3$. The resulting random tessellation of $W\subset{\Bbb R}^d$ is denoted by $Y_{\Lambda}(t,W)$ and is usually called the {\bf STIT tessellation} (for {\bf ST}able with respect to
 {\bf IT}erations) with driving measure $\Lambda$ and time threshold $t.$ By noting that for any for compact
 convex $V \subset W$ the trace the MNW-construction leaves on $V \subset W$ coincides in law with the MNW-process for $Y_{\Lambda}(t,V)$, we find that $Y(t,\cdot)$ enjoys a spatial consistency property in that $Y_{\Lambda}(t,W) \cap V \overset{D}{=} Y_{\Lambda}(t,V)$ , where $\overset{D}{=}$ stands for equality in distribution. Note that by facets, cells and regions of the restricted tessellation $Y_{\Lambda}(t,W) \cap V$ we understand the non-empty
 intersections of, respectively, facets, cells and regions of $Y_{\Lambda}(t,W)$ with $V.$ 
 The consistency readily allows us, by utilizing the consistency theorem \cite[Thm. 2.3.1]{SW},
 to construct the whole-space process $Y_{\Lambda}(t) := Y_{\Lambda}(t,{\Bbb R}^d)$ in such a way 
 that $Y_{\Lambda}(t,W) \overset{D}{=} Y_{\Lambda}(t) \cap W.$ Note that, by construction, all cells of $Y_\Lambda(t,W)$ and its global variant $Y_{\Lambda}(t)$ are necessarily convex (and even polyhedral if we assume the initial window $W$ to be a polytope). For terminological convenience we shall use the name {\bf region} for a general cell of the tessellation.\\ To justify the name {\it iteration stable} for the tessellation $Y_{\Lambda}(t)$ we note that if we insert an independent copy of the tessellation $Y_{\Lambda}(s) \cap c$
 into each active cell $c$ of the tessellation $Y_{\Lambda}(t),$ the resulting
 {\bf iterated tessellation} $Y_\Lambda(t) \boxplus Y_\Lambda(s),$ arising as the local superposition
 of i.i.d. copies of $Y_\Lambda(s)$ within the cells of $Y_\Lambda(t),$ coincides in law with $Y_\Lambda(t+s).$ 
 On the other hand, it is easily concluded from the MNW-construction that 
 \begin{equation}\label{SCALEINV}
  Y_{\Lambda}(t_1) \overset{D}{=} \frac{t_2}{t_1} Y_{\Lambda}(t_2),
 \end{equation}  
 which is to say $Y_{\Lambda}(t_1)$ coincides with $Y_{\Lambda}(t_2)$ re-scaled by the factor
 $t_2/t_1.$ Consequently,
 \begin{equation}\label{STABILITY}
  Y_{\Lambda}(t) \overset{D}{=} m Y_\Lambda(tm) \overset{D}{=} m(\underbrace{Y_\Lambda(t)\boxplus\ldots\boxplus Y_\Lambda(t))}_{m} \overset{D}{=}: m Y_{\Lambda}(t)^{\boxplus m}
 \end{equation}
 for $t > 0$ and $m \in {\Bbb N}.$ The relation (\ref{STABILITY}) is what usually
 goes under the name of stability with respect to iterations, whence the
 name STIT.\\ The stability relation with respect to iteration is one of the crucial properties for studying the
 MNW-construction. In fact it turns out \cite[Corollary 2]{NW05} that each
 translation invariant random tessellation of positive and finite surface density arising
 as a re-scaling of its iterations does admit a MNW-representation as above.
 Moreover, these are the only positive and finite surface density tessellations arising
 as re-scaled limits of iterations of stationary tessellations as shown by Theorem 3 in  
 \cite{NW05}. We like to point out that the STIT model enjoys truly remarkable properties with a lot of non-trivial, yet explicit, results available by now and under active development going on, compare with \cite{SLDP, ST, ST2, ST3, TW, NOW}.\\ \\ It is convenient from the formal viewpoint to regard the random STIT
 tessellation $Y_{\Lambda} = Y_{\Lambda}(t)$ as the random collection of maximal polytopes marked by
 their respective birth times, that is to say a marked point process in the space of
 $(d-1)$-dimensional compact polytopes in ${\Bbb R}^d$ carrying marks in $[0,t].$
 With probability one this collection
 \begin{description}
 \item{\bf[Empty start]} is empty at the time $0$, 
 \item{\bf [Local finiteness]} is locally finite in a sense that a compact subset of ${\Bbb R}^d$ intersects only a finite number of maximal polytopes,
  \item{\bf [Iterative binary subdivisions]} and enjoys the property that each facet (this is a maximal polytope) born at time $0 \leq s \leq t$ is contained in a unique convex cell and splits it into two convex sub-cells.
 \end{description}
 Each collection of time-marked $(d-1)$-dimensional polytopes with these three properties is what we call a
 {\bf nested tessellation} in this paper. Note, that unbounded cells and regions are permitted at this stage.
 The space of all nested tessellations of a given convex set $W \subseteq {\Bbb R}^d,$
 possibly with $W = {\Bbb R}^d,$ is denoted by ${\Bbb Y}_{[0,t]}(W)$
 and we abbreviate ${\Bbb Y}_{[0,t]} := {\Bbb Y}_{[0,t]}({\Bbb R}^d).$ 
 Endowing the space of maximal polytopes with the usual Hausdorff distance (see \cite[Chap. 12.3]{SW})
 allows us to define a natural vague topology on ${\Bbb Y}_{[0,t]}(W),$ see
 \cite[Definition A2.3.I]{DVJ88}, which induces the corresponding Borel $\sigma$-field on ${\Bbb Y}_{[0,t]}(W)$.
 Note that upon {\it forgetting} the time marks we may also interpret a tessellation
 $y \in {\Bbb Y}_{[0,t]}(W)$ as closed subset of ${\Bbb R}^d$ arising as the
 union of all its maximal polytopes. It is easily verified (compare with \cite[Thm. 12.2.3]{SW}) that this yields a continuous map from ${\Bbb Y}_{[0,t]}(W)$ to the space ${\cal F}(W)$ of closed subsets of $W$ with
 the usual Fell (hit-or-miss) topology \cite[Def. 1.1.1]{MOL05}. Another important
 measurable map is the time restriction map 
 $\iota_s : {\Bbb Y}_{[0,t]}(W) \to {\Bbb Y}_{[0,s]}(W),\; 0 \leq s \leq t,$ removing
 from its argument tessellation all maximal polytopes with birth times exceeding $s$ and keeping
 only those born no later than at the time $s.$ Note, that this, yet measurable, mapping is in general not continuous, which may be seen from the remarks on page 567 in \cite{SW}. By a \textbf{random nested tessellation} in $W \subseteq {\Bbb R}^d$ we shall obviously mean a
 random element taking values in the space ${\Bbb Y}_{[0,t]}(W)$ for some $t > 0.$ For a random nested
 tessellation $Y$ (note that we use $y$ for generic deterministic nested tessellations and
 $Y$ for generic random ones) we write 
 $$ Y(s) := \iota_s(Y),\; 0 \leq s \leq t, $$
 which yields tessellation-valued right-continuous stochastic process
 $(Y(s))_{s \in [0,t]}$ on ${\Bbb Y}_{[0,t]}(W)$ or ${\Bbb Y}_{[0.t]}$. We also write
 $$ Y(s,W) = Y(s) \cap W $$ for $s \in [0,t]$ and $W \subseteq {\Bbb R}^d.$ 
 Clearly, $(Y(s))_{s \in [0,t]}$ is a (possibly non-homogeneous) Markov
 process, because $Y(s)$ uniquely determines $Y(s')$ for all $s' < s.$ In this paper we will identify a nested random tessellion $Y\in{\Bbb Y}_{[0,t]}$ with the process $(Y(s))_{s\in[0,t]}$ induced by the time restriction map.\\ \\ To get more insight into the structure of these processes, we restrict our attention to random tessellations with laws locally absolutely
 continuous with respect to that of $Y_{\Lambda_{\rm iso}}(t),$ that is
 to say ${\cal L}(Y(t,W)) \ll {\cal L}(Y_{\Lambda_{\rm iso}}(t,W))$
 for all compact convex windows $W \subset {\Bbb R}^d$ (read ${\cal L}(...)$ as \textit{the law of ...}). Note in particular that this condition is satisfied for $Y = Y_{\Lambda}$ iff $\Lambda \ll \Lambda_{\rm iso}.$  For a fixed compact convex
 window $W \subset {\Bbb R}^d$ and $s \in [0,t],$ by the local absolute continuity
 as assumed above, the general theory of pure jump Markov processes \cite{EK}, \cite{LIG} guarantees that
 there exists a {\bf split kernel} $\Phi^W_s$ such that for each $y \in {\Bbb Y}_{[0,s]}(W)$
 and $c \in \Cells(y),$ $\Phi^W_s(\cdot|c,y)$ is a finite measure on $[c],$
 moreover $\Phi^W_s$ depends measurably on $c$, $y$ and $s$ and,
 furthermore, $\Phi^W_s$ is such that the process $(Y(s,W) = Y(s) \cap W)_{s \in [0,t]}$
 evolves in law according to the following recursive split dynamics
 \begin{description}
  \item{\bf [Recursive split dynamics]}
 \begin{itemize}
  \item Begin at the time $0$ with $Y(0,W) := \emptyset.$
  \item In time $s \in [0,t],$ for each $c \in \Cells(Y(s,W))$, with intensity 
        \begin{equation}\label{SplitInt}
         \Phi^W_s(dH|c,Y(s)) ds,
        \end{equation}
        split $c$ with hyperplane $H \in [c]$ selected according to the law \begin{equation}\Phi_s^W(dH|c,Y(s))\over |\Phi_s^W([c]|c,Y(s))|\label{selectH}\end{equation} into two new cells $c^+(H)$ and
        $c^-(H)$ (for definiteness say with $c^+(H)$ lying in the left half-space
        determined by $H$), thus adding to $Y(s,W)$ a new maximal polytope $H \cap c$ marked
        with the birth time $s.$ The resulting new tessellation is denoted
        by $y[c \oslash_s H].$ 
 \end{itemize}
\end{description}
 More formally, it means that $Y(s,W)$ is a (non-homogeneous) Markov process with generator ${\Bbb L}_s^W$ given by 
 \begin{equation}\label{LocGen}
  [{\Bbb L}^W_s f](y) = 
  \sum_{c \in \Cells(y)} \int_{[c]}
  [f(y[c \oslash_s H]) - f(y)]
  \Phi^W_s(dH|c,y),
 \end{equation}
 defined for all bounded measurable $f : {\Bbb Y}_{[0,t]}(W) \to {\Bbb R}.$  
 In particular, one can directly conclude from (\ref{selectH}) in view of (\ref{SplitMeasureSTIT}) that for the STIT tessellation $Y_{\Lambda}(t)$ we have 
 \begin{equation}\label{DlaPhiLambda}
  \Phi^W_{s;\Lambda}(dH|c,y) = {\bf 1}_{H \in [c]} \Lambda(dH),
 \end{equation}
 whence the kernel does not depend on $W$ for $c \subseteq W$ and also not on the time parameter $s$ and on the environment $y$.\\ In intuitive terms, all this means that a generic finite volume nested tessellation process 
 $Y(s,W),\; s \in [0,t],$ admits a representation where each cell splits with rate given by
 the total mass of its split kernel as depending on the surrounding configuration, whereas
 the probability law -- arising upon normalising this kernel -- governs the splitting mechanism, as 
 determining the cutting hyperplane. Importantly, the same informal description is also
 valid in infinite volume.
 \begin{proposition}\label{PHIINF}[Recursive split dynamics]
  There exist split kernels $\Phi_s$ such that for each $y \in {\Bbb Y}_{[0,s]}$
  and $c \in \Cells(y),$ $\Phi_s(\cdot|c,y)$ is a finite measure on the space $[c],$
  moreover $\Phi_s$ depends measurably on $c$ and $y$
  and, furthermore, $\Phi_s$ is such that the infinite-volume Markov
  process $(Y(s))_{s \in [0,t]}$ admits generator 
  \begin{equation}\label{OgGen}
   [{\Bbb L}_s f](y) = \sum_{c \in \Cells(y)} \int_{[c]}
   [f(y[c \oslash_s H]) - f(y)] \Phi_s(dH|c,y)
 \end{equation}
 defined for all bounded $f$ which are $\sigma(Y(s):s\in[0,t])$-measurable for some bounded 
 convex $W \subset {\Bbb R}^d.$ Moreover, for $H \in [c],$ 
 \begin{equation}\label{ZalMdzygen}
  \Phi^W_s(dH|c,y) = {\Bbb E}[\Phi_s(dH|Y(s))|Y(s) \cap W = y].
 \end{equation}
\end{proposition}
 Note that the formal statement (\ref{OgGen}) simply means that, in analogy to
 the finite volume set-up, also the infinite volume process $Y(s),\; s \in [0,t]$ unfolds in
 time according to the {\bf [Recursive split dynamics]}, where each cell $c$ splits by a hyperplane $H$
 with instantaneous intensity $\Phi_s(dH|c,y) ds$ at time $s$ in environment
 $y.$ The relation (\ref{ZalMdzygen}) states the intuitively obvious fact that for a cell $c$
 its instantaneous split rate by $H \in [c]$ coincides with the whole space split rate by $H,$
 averaged over the environment outside $W.$\\ A proof of Proposition \ref{PHIINF} in an extended setting is given in \cite{SLDP}, but we also give an argument in Appendix \ref{SecAp} in order to keep the paper self-contained.\\ \\ Whereas it is relatively easy to show that, as stated in Proposition \ref{PHIINF}, each nested tessellation admits a unique family of characterizing split kernels,
 the converse question whether a given family of kernels determines a tessellation
 and if this tessellation is unique and well behaved, is in general very difficult. In the
 present paper we shall handle this problem for a special class of so-called {\bf shape-driven Markov kernels}, which are introduced in the next section.

\section{Shape-Driven Markov Tessellations}\label{secMarkovTess}

\subsection{The Basic Setting}

 Proposition \ref{PHIINF} provides a very general structural characterization of random
 nested tessellations, at least those locally absolutely continuous with respect to the
 law of the STIT tessellation $Y_{\Lambda_{\rm iso}}(t).$ In 
 particular, it indicates -- in agreement with ideas recently developed by Cowan \cite{Cowan} --
 that it is natural to consider the splitting mechanism in iterative nested tessellations as 
 consisting of two separate ingredients
 \begin{itemize}
  \item \textit{split intensities}, corresponding to the total masses $|\Phi_s([c]|c,y)|$ of the split kernels,
  \item \textit{split geometry}, corresponding to the normalized kernels $\Phi_s(\cdot|c,y)/|\Phi_s([c]|c,y)|$ in our setting.
 \end{itemize}
 Whereas most of the work so far was concentrated on the first intensity-related point, while
 usually keeping the geometric ingredient of the kernels just homogeneous, compare with the so-called geometry
 independent apportionment considered in \cite{Cowan}, here our interest focuses on the second point.
 Our theory concentrates on a certain subfamily of admissible kernels chosen in such a way that the
 resulting class of tessellations is rich and flexible enough, but yet mathematically tractable. We begin
 with a natural spatial Markovianity postulate for the nested tessellation $Y$:
 \begin{description}
  \item {\bf [Markov property]} Given any two regions $A_1$ and $A_2$ of $Y(s)$, $0\leq s\leq t$, with disjoint interiors we have $Y(s,A_1)$ and $Y(s,A_2)$ conditionally independent.
 \end{description}
 This condition simply means that the behavior of the process within a cell -- and
 recursively, its sub-cells should they arise -- does depend on the cell only through its
 geometry and does not depend on what happens outside of it. In terms of split kernels
 of Proposition \ref{PHIINF} this is clearly equivalent to having the kernel $\Phi_s(dH|c,y)$
 depending only on $dH$ and $c$ but not $y$. That is to say, $\Phi_s(dH|c,y) = \Phi_s(dH|c)$
 for all environments $y,$ which even better justifies the use of notion of {\it Markovianity} in
 this context. \\ \\ Once we assume the so-understood spatially Markovian tessellations, it is important to
 note that if we just change the split intensities of the process while keeping its split
 geometry unchanged, that is to say if we multiply $\Phi_s(\cdot|\cdot)$ by scalar values,
 then we end up with another tessellation process amenable to the original, unmodified one
 by local time re-parametrization over individual regions, yet with the overall tessellation
 geometry remaining the same upon forgetting the time marks. This motivates us to
 stick to the usual {\it canonical} split intensities in this work, which enjoy particularly good
 properties due to their natural relationship to the previously developed STIT theory. Namely, in the sequel
 we shall always assume that the split intensity is given by the $\Lambda$-mean width
 (this is the usual mean width from convex geometry if $\Lambda = \Lambda_{\rm iso}$) as in the standard STIT setting:
 \begin{description}
  \item {\bf [Canonical split intensities]}
   We have 
   \begin{equation}\label{CSI}
    |\Phi_s([c]|c)| = \Lambda([c]).
  \end{equation}
 \end{description}
 One further requirement reducing unnecessary generality is to drop the dependency
 of $\Phi_s$ on $s,$ which is referred to as {\bf time homogeneity} with no need of
 further comments below. Finally, we want the split kernels $\Phi(\cdot|\cdot)$ to
 be isometry invariant and spatially homogeneous, where by the latter we mean: 
 \begin{description}
  \item{\bf [Spatial homogeneity]}
   We have
   \begin{equation}\label{SPH}
    \Phi(\alpha dH| \alpha c) = \alpha \Phi(dH|c),\;\; \alpha > 0.
   \end{equation}
 \end{description} 
 We shall say that a family of split kernels with canonical split intensities, which are isometry invariant
 and spatially homogeneous is {\bf shape-driven} and we shall stick to this level of generality
 in our considerations below. Note that the concept shape-driven as we introduced 
 above is meant to reflect the fact that the split mechanism considered here only takes into
 account the splitting cell's geometry, while setting apart its size, as opposed to the work done so far
 in the literature (cf. \cite{Cowan}), where the focus was mainly put on size and split intensities
 whereas the shape-related aspects of splits were often reduced to simple geometry-independent apportionments, see ibidem and Subsection \ref{SubsExa} below. A nested random tessellation constructed by a family of shape-driven split kernls having additionally the spatial Markov property is what we call a \textbf{shape-driven nested Markov tessellation} in this paper.

\subsection{Kernel Regularity and Main Results}
 To proceed, consider the natural space ${\cal K}$ of compact convex bodies in ${\Bbb R}^d$
 and let ${\cal K}_0 \subseteq {\cal K}$ be the subspace thereof, containing bodies {\it centered}
 at the origin $\bf 0,$ ${\cal K}_0 := \{ K \in {\cal K},\; c(K) = {\bf 0} \},$ where the {\it center} of
 a convex body $K$ stands here and below for some henceforth fixed shift-covariant
 selector of $K,$ say the Steiner point, the circumcenter, barycenter etc. Finally, we let
 ${\cal K}_{0;1} := \{ K \in {\cal K}_0,\; \Lambda([K]) = 1 \}$ be the family
 of centered convex bodies in ${\Bbb R}^d$ with unit split intensity,
 see (\ref{CSI}), also recall from Crofton's formula \cite[Thm 5.1.1.]{SW} that
 $\Lambda_{\rm iso}([K]) = \gamma_1 V_1(K)$ with $V_1$ standing for the order
 one intrinsic volume and with the integral-geometric constant $\gamma_1$ given by $$\gamma_1 = {2\kappa_{d-1}\over d\kappa_d} =\frac{\Gamma(d/2)}{\sqrt{\pi} \Gamma((d+1) / 2)},$$ where $\kappa_j$ stands for the volume of the $j$-dimensional unit ball. The split kernel $\Phi$ induces a natural Markovian (i.e. probabilistic)
 transition kernel $\Phi^{\oslash}$ (associated forward kernel, shrink kernel)
 from ${\cal K}$ into itself, acting by splitting a given $K \in {\cal K}$ into two cells
 according to the normalized kernel $\Phi(\cdot|\cdot)/|\Phi(\cdot|\cdot)|$
 and then outputting one of these cells chosen equiprobably at random.
 It is also convenient to consider the obvious modification $\Phi^{\oslash}_{0,1}$
 (re-normalized shifted associated forward kernel, re-normalized shifted shrink kernel)
 of $\Phi^{\oslash}$ acting from ${\cal K}_{0;1}$ into itself and
 obtained by suitably re-sizing and shifting the outcome of $\Phi^{\oslash}$ in
 the unique way to make it remain in ${\cal K}_{0;1}$. The obvious intermediate option,
 shifted shrink  kernel $\Phi^{\oslash}_0$ involving shifting but not re-scaling, shall also
 be considered in the sequel. The probabilistic kernels $\Phi^{\oslash}, \Phi^{\oslash}_0$
 and $\Phi_{0;1}^{\oslash}$ generate natural {\it shrink}, {\it shifted shrink} and
 {\it re-normalized shrink} Markov processes on ${\cal K}, {\cal K}_0$
 and ${\cal K}_{0;1}$, respectively, and we denote them by
 $$K^{\oslash}_1,K^{\oslash}_2,\ldots;\ \ \ \ K^{\oslash;0}_1, K^{\oslash;0}_2,\ldots\ \ \ \ \text{and}\ \ \ \ K^{\oslash;0;1}_1, K^{\oslash;0;1}_2,\ldots$$ with the law of $K^{\oslash}_{i+1}$
 given by $\Phi^{\oslash}(\cdot|K^{\oslash}_i),$ with $K^{\oslash;0}_{i}$ arising 
 as $K^{\oslash}_i$ shifted to the origin and, eventually, 
 with $K^{\oslash;0;1}_i$ arising as $K^{\oslash}_i$ re-sized and shifted to stay in
 ${\cal K}_{0;1}.$ Note that a shrink process should be regarded as 
 a realization in law of a randomly picked branch/chain in the successive
 forward cell division procedure of Proposition \ref{PHIINF} for the split kernel $\Phi.$
 Below it will be also convenient to consider the corresponding time-reversed growth
 process modeling the backwards history restoration along a random split branch,
 from tessellation cells back to their parents in successive sub-divisions.\\ 
 A continuous-time version of the shrink dynamics will also be in focus of our interest,
 namely we shall let $\hat{K}_t,\; t \geq 0$ be the ${\cal K}_0$-valued process 
 unfolding according to the following stationary Markovian dynamics in time $t$:
 \begin{description}
 \item{\bf [Continuous shrink dynamics]}
 \begin{description}
  \item{\bf [CSD1]}
          During the time interval $[t,t+dt]$ linearly re-size the body $\hat{K}^{\oslash;0}_t$
          with scaling factor $1+dt$, by keeping it in addition 
          centered at $\bf 0$ to stay within ${\cal K}_0.$
  \item{\bf [CSD2]} With intensity $\Lambda([\hat{K}^{\oslash;0}_t]) dt$ do split and
   replace $\hat{K}^{\oslash;0}_t$ according to the kernel $\Phi^{\oslash}_0.$
 \end{description} 
\end{description}
 The so-defined dynamics admits a number of stationary regimes which are, nevertheless,
 unique up to scaling under rather general regularity conditions, as we will see in the sequel.\\ Proceeding toward the formulation of the core of our theory, throughout this paper we
 shall impose some regularity conditions on $\Phi$ and $\Phi^{\oslash}$. We have decided to make them possibly somewhat stronger than needed, in order to keep the presentation and the argument below
 considerably simpler and more elegant without reducing the class of interesting
 examples and applications.
 \begin{description}
  \item{\bf [Split kernel regularity]}
  \begin{description} 
  \item{\bf [SKR1]} Let for any convex polytope $c\subset{\Bbb R}^d$, $f_c$ be a continuous density on $[c]$ with respect to some hyperplane measure $\Lambda\ll\Lambda_{\text{iso}}$ satisfying $f_{\alpha c}(\alpha dH)=f_c(dH)$ for any $\alpha>0$ and $f_{\phi(c)}=f_c$ for any isometry $\phi:{\Bbb R}^d\rightarrow{\Bbb R}^d$. We assume that the splitting hyperplane is chosen from the normalized hyperplane measure having this density, this is to say, the split kernels are of the form $$\Phi(dH|c)=f_c(H)\Lambda(dH).$$
  \item{\bf [SKR2]} The probabilistic transition kernel $\Phi^{\oslash}_{0;1}$ admits the unique invariant measure 
   $\varpi_{\Phi}$ on ${\cal K}_{0;1}$ for the stationary re-normalized shrink process $(K^{\oslash;0;1}_i)_{i \geq 1}.$
  \end{description}
 \end{description}
\begin{remark}
We expect that {\bf [SKR1]} already implies {\bf [SKR1]}, but we were not able to give a rigorous proof of this conjecture. 
\end{remark}
Recall that we say that a nested tessellation $(Y(s))_{s \in [0,t]}$ is compatible with a split kernel $\Phi$ iff the statement of Proposition \ref{PHIINF} holds with $\Phi_s = \Phi$ there. The main result of this paper is
  \begin{theorem}\label{MainThm}[Existence and uniqueness theorem]
   Assume that $\Phi$ is a Markovian shape-driven split kernel satisfying the split kernel regularity conditions from above. Then for each $t > 0$ there exists a unique nested tessellation $Y^\Phi=(Y^{\Phi}(s))_{s \in [0,t]}$ in ${\Bbb R}^d$ compatible with $\Phi.$ The tessellation $Y^{\Phi}$ enjoys the following properties:
   \begin{itemize}
    \item[(a)] $Y^{\Phi}(t)$ is spatially translation invariant (i.e. stationary).
    \item[(b)] The tessellation satisfies the scaling relation $Y^{\Phi}(t) \overset{D}{=} \alpha Y^{\Phi}(\alpha t),\; \alpha > 0$.
    \item[(c)] We have $Y^{\Phi}(s) = \iota_s(Y^{\Phi}(t))$ for $t>s,$ that is to say $Y^{\Phi}$ is consistent in time.
    \item[(d)] The typical cell $\TypicalCell(Y^\Phi(t))$ of $Y^{\Phi}(t)$ with distribution ${\Bbb Q}^\Phi(t)$ has finite mean volume and satisfies the scaling relation
                   \begin{equation}\label{TypDet}
                     {\Bbb Q}^{\Phi}(t) = [t^{-1}] \odot {\Bbb Q}^{\Phi},
                   \end{equation}
                   where $[\alpha] \odot$ stands for the usual spatial
                   scaling operation on the argument distribution with a positive factor $\alpha$ and
                   where ${\Bbb Q}^{\Phi}={\Bbb Q}^{\Phi}(1)$. 
                   Moreover, we have  
                   \begin{equation}{\Bbb E}\Lambda([\TypicalCell(Y^\Phi(t))])=
                   {d\over t}\label{FstOrd}\end{equation}
                   and, letting $\tilde{\Bbb Q}^{\Phi}(t)$ stand for the image of
                  the typical cell law ${\Bbb Q}^{\Phi}(t)$ under the standard re-normalizing
                  map ${\cal K}_0 \ni c 
                  \mapsto \tilde{c} = \frac{1}{\Lambda([c])} c \in {\cal K}_{0,1},$ 
                   \begin{equation}\label{TypCellRel1}
                     \tilde{{\Bbb Q}}^{\Phi}(1) = \tilde{{\Bbb Q}}^{\Phi} = \varpi_{\Phi}.
                   \end{equation}
                   Furthermore, the integral representation
            \begin{equation}\label{EqDistr}
             {\Bbb Q}^{\Phi}(t) = \int_{{\cal K}_0} ([t^{-1} \Lambda([c])] \odot \varpi_{\Phi})
             {\Bbb Q}^{\Phi}(dc).
            \end{equation}
            holds.
    \item[(e)] The typical cell distribution at the time $1,$ ${\Bbb Q}^{\Phi} = {\Bbb Q}^{\Phi}(1),$
                   is the unique stationary law for the {\bf [Continuous Shrink Dynamics]} satisfying
                   (\ref{FstOrd}) with $t=1$ there.
    \item[(f)] The the law ${\Bbb Z}^{\Phi}$ of the zero cell of $Y^{\Phi}(t)$ arises as the area-weighted
            modification of that of the typical cell, that is to say for $c \ni {\bf 0}$, 
            \begin{equation}\label{ZCell}
             {\Bbb Z}^{\Phi}(dc) = \frac{\Vol_d(c) {\Bbb Q}^{\Phi}(dc)}{\int_{{\cal K}_0}
                \Vol_d(c') {\Bbb Q}^{\Phi}(dc')}.
            \end{equation}
  \end{itemize}
 \end{theorem}
 As already mentioned before, the crucial feature of Theorem \ref{MainThm} is that the
 unique whole-space tessellation $Y^{\Phi}(t)$ gets constructed in a global construction
 greatly reminiscent of the (by now) well known Mecke-Nagel-Weiss (MNW) global construction
 developed in \cite{MNW} for the STIT processes. Hence, the connection between the
 incremental MNW-construction, represented here by the {\bf [Recursive split dynamics]},
 and the corresponding {\it global} MNW-construction, represented here by the global construction
 in our Subsection \ref{GlConstr}, carries over to the general class of shape-driven tessellations.
 It is also important to note that in view of our assumption on the stationary regimes for the continuous
 shrink dynamics, the point (e) of Theorem \ref{MainThm} combined with (\ref{FstOrd})
 allows to determine the typical cell ${\Bbb Q}^{\Phi}(t)$ and Markov-chain Monte-Carlo
 simulation can be an effective option.\\
 \begin{remark}\label{RemMeaningMainThm}
  It should be emphasized that, whereas our main Theorem \ref{MainThm} does guarantee
  the existence and uniqueness of the $\Phi$-compatible whole-space tessellation, it should not be
  mistakenly believed to give any kind of spatial consistency as in the STIT case. In general we should not expect
  that $Y^{\Phi}(t)\cap W$ coincide in law with the finite volume process we 
  would obtain if we carried out the recursive split construction restricted to a 
  proper subset $W \subset {\Bbb R}^d$ from the very beginning. Thus, unlike many other
  natural properties showing up in our general Theorem \ref{MainThm}, the consistency seems
  restricted for STIT and other closely related processes.
\end{remark}

\subsection{Examples}\label{SubsExa}
 We shall present now some classes of examples to which our theory may be applied:
\begin{enumerate}
 \item[(1)] We may consider the smooth constant density $\alpha\cdot{\bf{1}}[H\in[c]]$ with $\alpha>0$. This leads in view of (\ref{DlaPhiLambda}) to a time-changed STIT tessellation $Y^\Phi(t)$. More precisely, we have by (\ref{SCALEINV}) the relation $Y^\Phi(t)={1\over\alpha}Y_\Lambda(t)$, where $Y_\Lambda(t)$ is the standard STIT tessellation with driving hyperplane measure $\Lambda$.
 \item[(2)] For a convex set $K\subset{\Bbb R}^d$ and $r>0$ we denote by $\ero(K,r)$ the (possibly empty) erosion of $K$ and a ball with radius $r$, this is the set of all interior points of $K$ that have distance at least $r$ to the boundary of $K$. The natural density ${\bf 1}[H\in\ero(c,r)]$ is unfortunately not smooth on $s$. However, for any $\varepsilon>0$ there exists by Urysohn's classical lemma a continuous function $\hat{f}_c$ with $f_c=1$ on $[\ero(c,r)]$ and $\hat{f}_c=0$ on $[c\setminus\ero(c,r-\varepsilon)]$. Using this density, the split kernel $$\Phi(dH|c):=\hat{f}_c(H)\Lambda(dH)$$ is an admissible split-kernel in the sense of our theory, as long as $\Lambda\ll\Lambda_{iso}$. The restriction to the sets $\ero(c,r)$ implies that the splitting hyperplanes have at least distance $r-\varepsilon$ to the boundary of their mother cell, which is some kind of hard-core condition added to the construction of STIT tessellations.
\end{enumerate}
An interesting class of examples can be extracted from \cite{Cowan}, although Cowan's approach is entirely different from ours, as he works always in finite volumes and with a fixed number of cells. Ibidem, a random division rule of a cell $c$ was called a \textit{geometry-independent apportionment} if it only depends on the volume $\Vol_d(c)$ of $c$. Whenever a cell $c$ splits it creates sub-cells $c^+$ and $c^-$ and the volume proportion $U={\Vol_d(c^+)\over\Vol_d(c)}$ has distribution function $G(u|v)$ with $v=\Vol_d(c)$. Assume that $G$ is symmetric, i.e. $\displaystyle G(u|v)=1-\lim_{x\rightarrow 1-u}G(x|v)$, which is to say that $U$ and $\Vol_d(c^-)\over\Vol_d(c)$ are identically distributed. For normalization assume further $G(0|v)=0$ and $\displaystyle\lim_{u\rightarrow 1}G(u|v)=1$. We consider now a cell splitting mechanism with the property that $U$ as defined above follows the distribution function $G$ and that the dividing hyperplane is chosen with respect to a hyperplane measure $\Lambda\ll\Lambda_{iso}$. This corresponds to a geometry-independent apportionment in the sense of \cite{Cowan} and Theorem \ref{MainThm} ensures the existence of a whole-space nested Markov tessellation following that splitting rule. Two particular examples are
\begin{enumerate}
 \item[(3)] $G(u|v)=u$, the uniform distribution on $[0,1]$,
 \item[(4)] $\displaystyle G(u|v)={\bf 1}[0\leq u\leq 1]{\Gamma(2a)\over\Gamma(a)^2}\int_0^ut^{a-1}(1-t)^{a-1}dt$, a symmetric Beta-distribution with parameter $a>0$. Note that this example shares some features with Poisson-Dirichlet partitions considered in \cite[Chap. 2.2]{Bertoin} in the framework of fragmentation theory.
\end{enumerate}
A combination of these and other splitting mechanisms lead in view of the existence result from above to a wide range of constrained nested Markov tessellations with rather homogeneous cell sizes as needed for applications to real world problems.

\section{Geometry of Shape-Driven Nested Markov Tessellations}\label{secGeometry}

After we have ensured the existence of whole-space shape-driven nested Markov tessellations by Theorem \ref{MainThm}, we turn now to some aspects of their geometry. One aim is to derive some general mean values and mean value formulas for this class of random tessellations and to make connection and to add a broad class of examples to the recently developed theory from \cite{WC}. The reason, why we mostly restrict ourself to the particular interesting cases $d=2$ and $d=3$ is, firstly, that with increasing space dimension, the complexity of geometric objects determined by the tessellation grows rapidly. Secondly, many of the mean values depend on parameters, which are determined by the precise cell splitting mechanism and the number of these parameters and their meaning becomes quite complex.\\ \\ Let $Y=Y_\Phi(t)$ be a stationary shape-driven nested Markov tessellation in ${\Bbb R}^d$ whose split kernel satisfies the split kernel regularity conditions and denote by $X$ a generic class of convex polytopes determined by $Y$, for example the class $V$ of vertices or the class $C$ of its cells.\\ Let $X_1$ and $X_2$ be any two classes of convex polytopes. We say that $x_1\in X_1$ is \textbf{adjacent} to $x_2\in X_2$ if either $x_1\subseteq x_2$ or $x_2\subseteq x_1$. For $x_1\in X_1$, the number of objects of class $X_2$ adjacent to $x_1$ is denoted by $m_{X_2}(x_1)$. Moreover, the mean value $\mu_{X_1X_2}$ is defined as $$\mu_{X_1X_2}:={\Bbb E}_{X_1}m_{X_2}(x_1)=\int m_{X_2}(x_1)P_{X_1}(x_1),$$ where $P_{X_1}$ is the \textbf{Palm distribution} with respect to $X_1$, see \cite{SW} for exact definitions. We can interpret $\mu_{X_1X_2}$ as the mean number of objects of class $X_2$ adjacent to the \textbf{typical} object from class $X_1$. Furthermore, we let $\lambda_X$ be the \textbf{intensity} of objects of class $X$, i.e. the mean number of class-$X$-objects per unit volume and denote by $\nu_j(X)$ the mean number of $j$-faces of the typical object from class $X$, where $j$ ranges from $0$ to the dimension of the typical $X$-object.

\paragraph{The typical cell in the general isotropic case.}

Note that part (d) of Theorem \ref{MainThm} implies that ${\Bbb E}\Lambda(\TypicalCell(Y^\Phi(t)))$ is the same for all admissible split kernels $\Phi$ and is especially the same as for the typical cell of a Poisson hyperplane tessellation having intensity measure $t\Lambda$. Unfortunately, this does in general not imply that the mean intrinsic volumes of $\TypicalCell(Y^\Phi(t))$ coincide with that of the associated Poisson typical cell. However, for the special choice $\Lambda=\Lambda_{\text{iso}}$ we have $\Lambda([K])={2\kappa_{d-1}\over d\kappa_d}V_1(K)$ for any compact convex $K\subset{\Bbb R}^d$, which implies
\begin{equation}{\Bbb E}V_1(\TypicalCell(Y^\Phi(t)))={d^2\kappa_d\over 2\kappa_{d-1}}{1\over t},\label{V1TypicalCell}\end{equation} in particular $\pi/t$ in the planar and $6/t$ in the spatial case $d=3$.

\paragraph{The planar case.} Most of the first-order geometry of shape-driven nested Markov tessellations in ${\Bbb R}^2$ can be determined by the vertex geometry of these tessellations. It follows directly from the construction and our assumption that the law of $Y^\Phi(t)$ is absolutely continuous with respect to that of the STIT tessellation $Y_{\Lambda_{\rm iso}}(t)$, that such tessellations can -- in the planar case -- only have T-shaped vertices, which is to say that from any vertex we have exactly $3$ outgoing edges and two of them are collinear. Here and in the sequel, we understand by an \textbf{edge} a line segment bounded by two vertices, but with no further vertices in its relative interior. In contrast, by a \textbf{side} of the tessellation, we understand any of the cell's sides. These line segments are of course bounded by vertices, but may have further vertices in their relative interiors. The classes of edges and sides of the tessellation are denoted by $E$ and $S$, respectively. Moreover, we consider the class $I$ of its I-segments introduced in Section \ref{secSTIT}. Observe now that any vertex of the tessellation is ...
\begin{itemize}
 \item[-] ... endpoint of exactly $3$ edges,
 \item[-] ... endpoint of exactly $4$ sides,
 \item[-] ... endpoint of exactly $1$ I-segment and
 \item[-] ... corner of exactly $2$ cells.
\end{itemize}
Using these properties, we can relate the intensities $\lambda_E$, $\lambda_S$, $\lambda_I$ and $\lambda_C$ to the vertex intensity $\lambda_V$ as follows: $$\lambda_E={3\over 2}\lambda_V,\ \ \ \lambda_S=2\lambda_V,\ \ \ \lambda_I={1\over 2}\lambda_V\ \ \ \text{and}\ \ \ \lambda_C={1\over 2}\lambda_V.$$ Moreover, denoting the mean lengths of the different typical line segments by $L_E$, $L_S$ and $L_I$, we get \begin{equation}L_I=2L_S=3L_E,\label{MeanValues2d}\end{equation} by noting that each interior point of a side is counted twice. Moreover, application of the standard mean value formulas from \cite{SW} or \cite{SKM} leads to $\mu_{CV}=2+2\lambda_V/\lambda_C=6$, i.e. the typical cell has $6$ vertices on its boundary in the mean. On the other hand, it is well known that the mean number of corners $\nu_0(C)$ of the typical cell satisfies $$\nu_0(C)=\nu_1(C)={2(\mu_{VE}-1)\over\mu_{VE}-2}={2(3-1)\over 3-2}=4,$$ which can be interpreted by saying that the typical cell is a quadrangle in the mean with two additional vertices on its boundary, which are no corners of it.\\ In the special isotropic case $\Lambda=\Lambda_{\text{iso}}$ we can combine (\ref{V1TypicalCell}) with the fact that the first intrinsic volume of a convex set equals half of its perimeter length, to obtain $$p=2{\Bbb E}V_1(\TypicalCell(Y^\Phi(t)))={2\pi\over t}$$ for the mean perimeter length $p$ of the typical cell of $Y^\Phi(t)$. Combined with $\nu_0(C)=4$ this leads to precise values for $L_I$, $L_S$, $L_E$ and to the observation that in the planar case the edge length density $L_A$ equals the construction time $t$ (this is no longer true for $d>2$). Standard mean value formulas from \cite{SKM} also imply exact values for the intensities $\lambda_I$, $\lambda_S$, $\lambda_S$, $\lambda_E$ and $\lambda_V$. Further, the mean area $a=1/\lambda_C$ of the typical cell can be found. The results are summarized in the following table
\begin{center}
\begin{tabular}{|c|c|c|c|c|c|c|c|c|c||c|c|c|}
\hline
$L_A$ & $\lambda_V$      & $\lambda_E$      & $\lambda_S$      & $\lambda_I$ & $L_E$ & $L_S$ & $L_I$ & $p$ & $a$ & $\mu_{VE}$ & $\mu_{CV}$ & $\nu_0(C)$\\
\hline
$t$   & ${2\over\pi}t^2$ & ${3\over\pi}t^2$ & ${4\over\pi}t^2$ & ${1\over\pi}t^2$ & ${\pi\over 3}{1\over t}$ & ${\pi\over 2}{1\over t}$ & $\pi{1\over t}$ & $2\pi{1\over t}$ & ${\pi\over 2}{1\over t^2}$ & $3$ & $6$ & $4$\\
\hline
\end{tabular}
\end{center}
Note, that these mean values coincide with those for stationary and isotropic STIT tessellations in the plane, see \cite{NW06}, which shows that the first-order geometry of isometry invariant planar shape-driven nested Markov tessellations is rather rigid due to their topological structure. However, we expect that the additionally introduced flexibility has an important influence on the second-moment structure of these tessellations.

\begin{figure}[t]
\begin{center}
\includegraphics[scale=0.6]{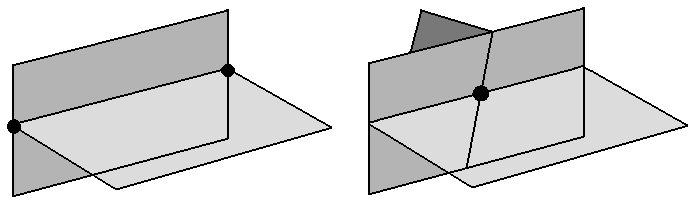}
\label{figvert}
\caption{The two possible types of vertices: T-vertices (left) and X-vertices (right)}
\end{center} 
\end{figure}

\paragraph{The spatial case.} According to the recent theory presented in \cite{WC}, a number of the most basic topological mean values of a stationary spatial random tessellation is determined by a system of $10$ parameters $\mu_{VE}$, $\mu_{EP}$, $\mu_{PV}$, $\mu_{CV}$, $\mu_{CE}$, $\mu_{CP}$, $\xi$, $\kappa$, $\psi$ and $\tau$, where $P$ stands for the class of \textbf{plates} of the tessellation and by a plate we mean a $2$-dimensional convex polygon bounded by edges but with no vertices or edges in its relative interior. Remarkably, this system of parameters can considerably be reduced in our case. To do so, we observe at first that the vertices of a shape-driven nested Markov tessellation in ${\Bbb R}^3$ can only be of two types: T-vertices and X-vertices, see Figure \ref{figvert}. This follows, similarly as in the planar case, from the fact that the law of $Y^\Phi(t)$ is absolutely continuous with respect to that of the STIT tessellation $Y_{\Lambda_{\rm iso}}(t)$.
Now, we can give an interpretation of the four different Greek parameters mentioned above. The parameter $\kappa$ is the proportion of T-vertices present in the tessellation, whereas $\xi=1$, because all edges are so-called $\pi$-edges in the sense of \cite{WC}. The parameter $\psi$ is the expected number of cell-side interiors adjacent to the typical vertex and $\tau$ is the expected number of plate-side interiors adjacent to the typical vertex. From the geometry of the vertices we can see that any T-vertex is adjacent to exactly one cell-side interior and to exactly one plate-side interior. Moreover, any X-vertex is adjacent to exactly $4$ cell-side interiors and to exactly $2$ plate-side interiors, so that we get $$\psi=1\cdot\kappa+4\cdot(1-\kappa)=4-3\kappa,\ \ \ \tau=1\cdot\kappa+3\cdot(1-\kappa)=1+\kappa.$$ Moreover, we obviously have $\mu_{VE}=4$ and  $\mu_{EP}=3$. Also the parameters $\mu_{CV}$, $\mu_{CE}$ and $\mu_{CP}$ can be expressed by $\kappa$ and $\mu_{PV}$ by using the standard mean value formulas for spatial random tessellations from \cite{SKM}. Thus, the only remaining parameters in our case are $\kappa$ and $\chi:=\mu_{PV}$.\\ Denote by $\lambda_{V[T]}$ and $\lambda_{V[X]}$ the intensity of T- and X-vertices, respectively, i.e. $\lambda_{V[T]}=\kappa\lambda_V$ and $\lambda_{V[X]}=(1-\kappa)\lambda_V$. We use now the vertex geometry to obtain the intensities $\lambda_I$, $\lambda_{S[C]}$, $\lambda_{S[P]}$ and $\lambda_E$ of I-segments, ${S[C]}$-segments, ${S[P]}$-segments and edges, respectively. Here, by an edge we mean again a line segment bounded by vertices with no vertices in its relative interior. A ${S[C]}$-segment is the ($1$-dimensional) side of a cell and a ${S[P]}$-segment is the side of a plate (thus, the classes $S[C]$ and $S[P]$ are multi-sets of line segments with multiplicities $1$ and $2$ and $1$, $2$ and $3$, respectively). Moreover, the I-segments are the sides of the maximal polygons inserted during the construction of the tessellation. We observe now that each T-vertex (X-vertex) is the endpoint of ...
\begin{itemize}
 \item[] ... exactly $2$ ($0$) I-segments,
 \item[] ... exactly $4$ ($4$) edges,
 \item[] ... exactly $6$ ($0$) $S[C]$-segments,
 \item[] ... exactly $10$ ($8$) $S[P]$-segments.
\end{itemize}
Using these facts we get: $\lambda_I={1\over 2}\cdot 2\lambda_{V[T]}=\kappa\lambda_V$, $\lambda_E={1\over 2}\cdot 4\lambda_V=2\lambda_V$, $\lambda_{S[C]}={1\over 2}(10\lambda_{V[T]}+8\lambda_{V[X]})=(4+\kappa)\lambda_V$ and $\lambda_{S[P]}={1\over 2}(6\lambda_{V[T]}+0\lambda_{V[X]})=3\kappa\lambda_V$. Furthermore, the mean lengths $L_I$, $L_{S[C]}$, $L_{S[P]}$ and $L_E$ can be calculated from $\lambda_IL_I=\lambda_EL_E=L_V$, $\lambda_{S[C]}L_{S[C]}=2L_V$ and $\lambda_{S[P]}L_{S[P]}=3L_V$, where $L_V$ is the edge-length density of the tessellation, so that we end up with the general relationship $$L_I={3\over 2}L_{S[C]}={4+\kappa\over 3\kappa}L_{S[P]}={2\over\kappa}L_E,$$ which in particular shows that the mean number of vertices in the relative interior of the typical I-segment equals ${2-\kappa\over\kappa}$, whereas $$\mu_{\text{relint}(S[C])V}={L_{S[C]}\over L_E}-1={4-3\kappa\over 3\kappa}\ \ \ \text{and}\ \ \ \mu_{\text{relint}(S[P])V}={L_{S[P]}\over L_E}-1={2-\kappa\over\kappa+4},$$ where $\text{relint}(X)$ stands for the relative interior of the typical object of class $X$.\\ The geometry of the classes of objects of dimension $2$ is considerably more involved and beside $\kappa$ also our second parameter $\chi$ enters the expressions. From the general formulas in \cite[Table 2]{WC} we obtain at first $\lambda_C={6-\chi\over\chi}\lambda_V$ and with the Euler-type relation $\lambda_V-\lambda_E+\lambda_P-\lambda_C=0$ it follows $\lambda_P={6\over\chi}\lambda_V$. Moreover, we have $$\lambda_F=\left({12\over\chi}+\kappa-2\right)\lambda_V.$$ Note, that by a facet we mean each facet of a cell, which implies that the class $F$ of all tessellation facets is a multi-set with multiplicities $1$ and $2$. This together with $\lambda_PA_P=S_V$ and $\lambda_FA_F=2S_V$ implies now $$A_F={12\over 12+\chi(\kappa-2)}A_P.$$ Moreover, we can calculate $\nu_0(P)=\nu_1(P)$ and $\nu_0(F)=\nu_1(F)$, the number of corners of the typical plate and that of the the  typical facet: $$\nu_0(P)=\nu_1(P)={\chi\over 6}(4+\kappa),\ \ \ \nu_0(F) = \nu_1(F) = {6\kappa\chi\over 12+\chi(\kappa-2)},$$ see \cite[Table 6]{WC}. Observe now that any edge of the tessellation is contained in the relative interior of precisely one cell facet, whence $$\mu_{\text{relint}(F)E}={\lambda_E\over\lambda_F}={2\chi\over 12+\chi(\kappa-2)}.$$ Next, we note that any T-vertex is contained in the relative interior of exactly one facet and that there is no facet in whose relative interior an X-vertex is contained, which implies $$\mu_{\text{relint}(F)V}={\kappa\chi\over 12+\chi(\kappa-2)}$$ and the general relationship ${\kappa}\mu_{\text{relint}(F)E}=2\mu_{\text{relint}(F)V}$. Moreover, $\mu_{FP}=2\lambda_P/\lambda_F=12/(12+\chi(\kappa-2))$, because each plate is part of two facets. Unfortunately and in contrast to the plates and facets, the mean values for maximal polygons cannot be expressed by the two parameters $\kappa$ and $\chi$ and we skip their discussion for this reason. For the approach in the STIT case we cite \cite{TW} and remark that the method developed there cannot be extended to any other shape-driven nested Markov tessellation.\\ To derive formulas for the typical cell $C$ of the tessellation $Y$, we can utilize the standard mean value formulas for spatial tessellations from \cite[Chap. 10.4]{SKM}, which lead to $$\mu_{CV}={4\chi\over 6-\chi},\ \ \ \ \mu_{CE}={6\chi\over 6-\chi},\ \ \ \ \mu_{CP}={12\over 6-\chi}$$ and furthermore the formulas from \cite[Table 6]{WC} to obtain $$\nu_0(C)={2\kappa\chi\over 6-\chi},\ \ \ \nu_1(C)={3\kappa\chi\over 6-\chi},\ \ \ \nu_2(C)={12+(\kappa-2)\chi\over 6-\chi}$$ for any shape-driven nested Markov tessellation in ${\Bbb R}^3$.\\ From our results, we can deduce some inequalities for the basic parameters $\kappa$ and $\chi$. First, observe that from $$\nu_0(P)={\chi\over 6}(4+\kappa)\geq {2\over 3}\chi\geq 3$$ we get $\chi\geq{9\over 2}$. Moreover, we obviously have $\chi<6$, because $\mu_{CP}$ for example has to be positive. Moreover, $\nu_0(C)\geq 4$ implies $$\kappa\geq{12-2\chi\over\chi}.$$ 
\begin{remark}
It would be interesting to know if $\chi$ can be arbitrarily close to $6$ and if there are examples for which $\kappa$ becomes arbitrarily close to $1$. These questions are closely related to open problems raised in \cite{WC}. 
\end{remark}

\paragraph{The STIT-case.} Mean values for planar and spatial random STIT tessellations were considered in \cite{NW06}, \cite{NW08} and \cite{TW}. Using the additional iteration stability and the property that STIT tessellations have Poisson-typical cells, the parameters $\lambda_{V}$ and $L_A$ in the planar case and $\lambda_V$, $\kappa$, $\chi$, $L_V$ and $S_V$ in the spatial case can be made available. In the planar case, the two parameters are $L_A=t$ and $\lambda_V=t^2\zeta$ with $\zeta$ given by $$\zeta=\int_{{\cal S}_1}\int_{{\cal S}_1}[u_1,u_2]{\cal R}(du_2){\cal R}(du_1),$$ where $\cal R$ is the directional distribution introduced in Section \ref{secSTIT} and where $[u_1,u_2]$ stands for the area of the parallelogram spanned by $u_1$ and $u_2$. In the spatial case, we derive from $\cal R$ the two constants
\begin{eqnarray}
\nonumber \zeta_2 &=& \int_{{\cal S}_2}\int_{{\cal S}_2}[u_1,u_2]{\cal R}(du_2){\cal R}(du_1),\\
\nonumber \zeta_3 &=& \int_{{\cal S}_2}\int_{{\cal S}_2}\int_{{\cal S}_2}[u_1,u_2,u_3]{\cal R}(du_3){\cal R}(du_2){\cal R}(du_1),
\end{eqnarray}
where, similarly to the planar case, $[u_1,u_2,u_3]$ is the volume of the parallelepiped spanned by $u_1$, $u_2$ and $u_3$. By a direct verification it can be shown that in the isotropic case, $\zeta_2$ and $\zeta_3$ are given by $\pi\over 4$ and $\pi\over 8$, respectively. We are now in the position to express the parameters discussed above by the construction time $t$ and the two numbers $\zeta_2$ and $\zeta_3$: $$S_V=t,\ \ L_V=S_V\zeta_2,\ \ \lambda_V=S_V^3\zeta_3,\ \  \kappa={2\over 3},\ \ \chi={36\over 7}.$$

\section{The Global Construction and Proof of Theorem \ref{MainThm}}\label{secProof}

 As already signaled above, our proof of the principal Theorem \ref{MainThm} goes by providing
 first a global construction for the infinite volume $\Phi$-splitting tessellation, constituting a direct
 and natural extension of the procedure given by Mecke, Nagel and Weiss in \cite{MNW, MNW2} for
 STIT tessellations. We then show that the so-obtained tessellation is the unique one compatible
 with $\Phi$ in the full space ${\Bbb R}^d$ and we verify a number of properties thereof, obtained
 as by-products of the argument.

 \subsection{The Global Construction}\label{GlConstr}

  The idea underlying our global construction involves a gradual reconstruction of the space-tessellating
  field, given its parts already reconstructed, starting from the void and eventually covering the 
  whole of ${\Bbb R}^d.$ To this end, we consider first the time-reversal of the shrink
  processes and the associated Markovian (probabilistic) growth kernels. For $k,k' \in {\cal K}_{0;1}$ we let
  \begin{equation}\label{OutgrowthRenormalized}
   \Phi^{*;0;1}(dk|k') = {\Bbb P}(K_i^{\oslash;0;1} \in dk|K_{i+1}^{\oslash;0;1}=k') =
   \frac{\Phi^{\oslash;0;1}(dk'|k) \varpi_{\Phi}(dk)}{\varpi_{\Phi}(dk')},
  \end{equation}
  which is defined $\varpi_{\Phi}$-almost everywhere and is referred to as the 
  {\bf re-normalized growth kernel} in the sequel. Dropping the re-normalization
  we also define, with $k,k' \in {\cal K},$
  \begin{equation}\label{OutgrowthAbsolute}
   \Phi^*(dk|k') = {\Bbb P}(K^{\oslash}_i \in dk|K_{i+1}^{\oslash}=k') =
   \frac{\Phi^{\oslash}(dk'|k) \varpi_{\Phi}(d \langle k \rangle_{0;1})}{\varpi_{\Phi}(d \langle 
   k' \rangle_{0;1})}
  \end{equation}
  where $\langle k \rangle_{0;1}$ stands for the standardized version of $k$,
  re-sized and shifted so as to fall into ${\cal K}_{0;1}.$  Below we grant the name
  {\bf growth kernel} to the so-defined $\Phi^*.$  Given the Markovian growth kernel
  $\Phi^*$ and a {\bf seed cell} $S_{\rm seed} \in {\cal K}$ we consider the
  {\bf spinal growth process} $S^*_0 = S_{\rm seed}, S^*_1, S^*_2,\ldots$, 
  corresponding in law to the time-reversed shrink process with infinite
  indexation, satisfying
  \begin{equation}\label{SpinalTransition}
   {\Bbb P}(S^*_{i+1} \in dk|S^*_i = k') = \Phi^*(dk|k').
  \end{equation}
 For the so-defined process, we clearly have by definition $S^*_{i+1} \supset S^*_i$ almost surely for $i=0,1,2,\ldots$.\\  To make proper use of the spinal growth process we shall now ascribe auxiliary {\it time marks} $S^*_i \mapsto \tau^*_i$  to its subsequent cells, with the objective of having 
  the joint process $(S^*_i,\tau^*_i)_{i \geq 0}$ described in law as follows:
  \begin{itemize}
   \item $S^*_i$ \--- the subsequent ancestor cells, backwards in time along a randomly
                           picked branch of the recursive split procedure generating the
                            tessellation considered,
   \item $\tau^*_i$ \--- the corresponding cell's birth times in the course of the recursive
                               split procedure, where by birth time $\tau^*_i$ of cell $S^*_i$ we
                               understand the time moment in $[0,t]$, where the super-cell
                               $S^*_{i+1}$ of $S^*_i$ splits and gives rise to $S^*_i.$ Clearly,
                               $0<\ldots\tau^*_{i+1} \leq \tau^*_i \ldots \leq t$ almost surely.
 \end{itemize}
  By convention, we put $\tau^*_0 = t$ for the seed cell. Then, aiming at
  the above interpretation and recalling (\ref{CSI}) determining the total split intensity
  $|\Phi([k_i]|k_i)| = \Lambda([k_i])$ for a cell $k_i$,
  we conclude that we have, conditionally on $S^*_i, S^*_{i+1}$ and
  $\tau^*_i$ given,
  $$ {\Bbb P}(\tau^*_{i+1} \in dt_{i+1}|S^*_{i+1}=k_{i+1}, S^*_i=k_i, \tau^*_i=t_i)
      \propto e^{-\Lambda([k_i])(t_i-t_{i+1})} dt_{i+1}, $$
  because the corresponding conditional probability of surviving the time period
  $[t_{i+1},t_i)$ in state $k_i$ upon shrinking to it at $t_{i+1}$ is just
  $\Lambda([k_i]) e^{-\Lambda([k_i])(t_i-t_{i+1})}dt_{i+1}$ by the construction introduced in Section \ref{secSTIT}. 
  Solving the proportion yields 
  \begin{equation}\label{CondTMark}
   {\Bbb P}(\tau^*_{i+1} \in dt_{i+1}|S^*_{i+1}=k_{i+1}, S^*_i=k_i, \tau^*_i=t_i) =
     \Lambda([k_i]) e^{-\Lambda([k_i])(t_i-t_{i+1})} dt_{i+1}.
  \end{equation} 
  Consequently, by (\ref{CondTMark}) above, the extended Markovian kernel
  \begin{equation}\label{JointTMGrow}
   \hat{\Phi}^*(dk,dt|k',t') = \Phi^*(dk|k') \Lambda([k'])\exp(-\Lambda([k'])(t'-t)) dt 
  \end{equation}
  governs the joint Markovian dynamics of the time-marked spinal growth process 
  $(S^*_i,\tau^*_i)_{i \geq 0}.$\\  We close this paragraph by noting that the cells $S_i^*$ eventually cover the whole space, which is $\bigcup_{i=0}^{\infty} S^*_i = {\Bbb R}^d$. This can be seen as follows, see also \cite[Lem. 4.1]{MNW}: Denote by $\sigma<t$ the largest birth time of a cell from the spinal growth process hitting the convex hull of $B_1(p)$ and the center $c(S_{\rm seed})$ of the seed cell, where $B_1(p)$ is the ball around $p\in{\Bbb R}^d$ with radius $1$. Then $p$ is contained in the interior of $S_{\tau_i^*}^*$ for any $\tau_i^*<\sigma$, thus the event $\{\exists i:B_1(p)\subset S_{\tau_i^*}^*\}$ has probability one. The claim follows now from the fact that ${\Bbb R}^d$ can be covered by balls $B_1(p)$ form a countable set of midpoints $p$ and the fact that the process $S_i^*$ is growing, that is $S_{i}^*\subset S_{i+1}^*$.

  \paragraph{The Construction.}
  The following crucial procedure is our {\bf global construction} with initial seed cell
  $S_{\rm seed},$ outputting the tessellation $Y^{\Phi}_{S_{\rm seed}}(t),\; t > 0.$ It can be regarded as a direct generalization of the global construction of STIT tessellations from \cite{MNW, MNW2}. 
  \begin{description}
   \item{\bf [Construction: Spinal phase]} Starting from the initial seed cell $S_{\rm seed}$ at
    terminal time $t$, construct the time-marked spinal growth process
    $(S^*_i,\tau^*_i)_{i \geq 0}$ backwards in time as determined by the extended spinal
    growth kernel $\hat{\Phi}^*$, see (\ref{JointTMGrow}). In this way a (non-stationary) {\bf frame
    tessellation} of the entire ${\Bbb R}^d$ is obtained, consisting of an infinite chain of nested cells
   fully covering the space.
  \item{\bf [Construction: Local fill phase]} Given the spinal chain with time-marked cells, for
          each cell $c$, having time mark falling below the threshold $t$, carry out within $c$
          the usual {\bf [Recursive split dynamics]} with split kernel $\Phi$ up to time $t.$ Output
          the resulting tessellation as $Y^{\Phi}_{S_{\rm seed}}(t).$
 \end{description}

 \subsection{Completion of the Proof}
  The further argument in our proof of Theorem \ref{MainThm} relies on showing that
  the desired unique tessellation $Y^{\Phi}$ can be obtained from the global construction from
  the last subsection under appropriate (randomized) choice of the initial seed cell $S_{\rm seed}$
  made coincide with the corresponding zero cell for $Y^{\Phi}.$  The reasoning
  splits into several steps and goes as follows.

  \paragraph{The existence of a translation invariant field.}
 We identify a random tessellation in ${\Bbb R}^d$ with the union of its cell boundaries, i.e. we regard it as a random element in the space $\cal F$ of closed subsets of ${\Bbb R}^d$. For our split kernel $\Phi$, satisfying our imposed regularity assumptions, we consider the family $${\cal Y}=\left\lbrace {\cal L}(Y^\Phi(t,[-R,R]^d)):\ R \in {\Bbb N} \right\rbrace $$ of laws of tessellations $Y^\Phi(t,[-R,R]),\;
R \in {\Bbb N}$. As the space $\cal F$ is compact, the familiy $\cal Y$ is tight. Whence, by Prohorov's theorem, there exists a convergent subsequence in $\cal Y$, whose limit is the distribution of a random closed set $Y^\Phi(t;\infty)$ in ${\Bbb R}^d$. We note that $Y^\Phi(t;\infty)$ is the frame of a random tessellation, the cell interiors of which are the connected components of $Y^\Phi(t;\infty)$'s complement. Moreover,
$Y^{\Phi}(t;\infty)$ is $\Phi$-compatible, as obviously the kernel $\Phi$ has the Feller property due to our assumption \textbf{[SKR1]} on the densities $f_c$. Consider next the family $$ {\cal Y}'  =\left\lbrace \frac{1}{(2R)^d} \int_{[-R,R]^d} [Y^\Phi(t,\infty) + x] dx  
    :\ R \in {\Bbb N} \right\rbrace $$
 of versions of $Y(t;\infty)$ smeared (shifted and normalized) over increasing cubes. Using
 again the tightness and passing to a subsequence we find in the limit $Y^{\Phi}(t)$, which is $\Phi$-compatible and translation invariant. Moreover, we readily conclude from the construction
  that $Y^{\Phi}$ can be taken consistent in time and that the desired standard scaling relation 
  \begin{equation}\label{STDscaling}
  \alpha Y^{\Phi}(\alpha t) \overset{D}{=} Y^{\Phi}(t),\; \alpha > 0,
  \end{equation}
  holds in view of (\ref{CSI}), possibly upon further passing to a subsequence.

 \paragraph{The typical cell and the zero cell.}
  Write $\lambda_C^{\Phi}(t)$ for the cell density of $Y^{\Phi}(t)$ and let ${\Bbb Q}^{\Phi}(t)$
  be the corresponding typical cell distribution on ${\cal K}_0$, 
  so that the standard cell intensity measure
  decomposition \cite[Eq. (4.2)]{SW} becomes
  $$ \Theta^{\Phi}(t) = \lambda_C^{\Phi}(t) {\rm Shift}(\ell^d \otimes {\Bbb Q}^{\Phi}(t)) $$
  where $\Theta^{\Phi}(t) = {\Bbb E} \sum_{c \in \Cells(Y^{\Phi}(t))} \delta_{c}$ is 
  the particle intensity measure in the sense of \cite[Section 4.1]{SW}, $\ell^d$ is the
  usual $d$-dimensional Lebesgue measure, whereas ${\rm Shift}(\cdot)$ stands for transporting
  the argument measure via the shift map $(x,c) \mapsto x+c$ with $x \in {\Bbb R}^d$ and $c \in
  {\cal K}_0.$  Under this notation, the {\bf [Recursive split dynamics]} of Proposition \ref{PHIINF}
  yields in view of (\ref{CSI})
  $$\frac{d}{dt}  \lambda_C^{\Phi}(t) = \lambda_C^{\Phi}(t) \int_{{\cal K}_0}
      \Lambda([c]) {\Bbb Q}^{\Phi}(t) (dc) = \lambda_C^{\Phi}(t){\Bbb E} \Lambda([\TypicalCell(Y^{\Phi}(t))]).$$ Indeed, this follows readily from the definition of the generator together with the definition of the typical cell, see \cite[Eq. (4.6) and Thm. 4.1.3]{SW}. However, by stationarity and by the scaling relation (\ref{STDscaling}), there exists a constant $a>0$ not depending on $t$ such that $\lambda_C^\Phi(t)=at^d$.
 Inserting this expression into the last equation immediately leads to \begin{equation}\int_{{\cal K}_0}\Lambda([c]){\Bbb Q}^\Phi(t)(dc)={\Bbb E}\Lambda([\TypicalCell(Y^{\Phi}(t))])={d\over t},\label{SplitBalanceI}\end{equation} which was part (\ref{FstOrd}) of statement (d) in Theorem \ref{MainThm}. In particular, this shows that the size of the typical cell of $Y^\Phi(t)$ -- measured here by $\Lambda([\cdot])$ -- is influenced only by the split intensities (\ref{SplitInt}), i.e. the total masses of the split kernel $\Phi$, and not by the precise split geometry determined by it. Moreover, it ensures in view of the well-known Usryson inequality from convex geometry (cf. \cite[Thm. 8.9]{Gruber}) that the typical cell has finite mean volume, adding to the statements of part (d) of the Theorem.\\ By the same argument as used above, again by (\ref{OgGen}) and (\ref{CSI}), we have
  \begin{eqnarray}  
  \nonumber \frac{d}{dt} \Theta^{\Phi}(t) &=&  \lambda_C^{\Phi}(t) {\rm Shift}\left( \ell^d \otimes
       \int_{{\cal K}_0}
       [2\Phi^{\oslash}(\cdot|c)-\delta_c] \Lambda([c]) {\Bbb Q}^{\Phi}(t)(dc) \right)\\
  \nonumber &=& \lambda_C^{\Phi}(t) {\rm Shift}\left( \ell^d \otimes 
      {\Bbb E}[2\Phi^{\oslash}(\cdot|\TypicalCell(Y^{\Phi}(t)))-
       \delta_{\TypicalCell(Y^{\Phi}(t))}] \right),
  \end{eqnarray}
  where the inner bracket reflects the operation of replacing the splitting cell $c$ by the
  two new cells arising. Putting the equations for $\lambda_C^\Phi(t)$ and $\Theta^\Phi(t)$ together yields 
  \begin{equation}\label{TypCellEqn}
   \frac{d}{dt} {\Bbb Q}^{\Phi}(t) = \int_{{\cal K}_0} [2\Phi^{\oslash}_0(\cdot|c) - \delta_c 
    - {\Bbb Q}^{\Phi}(t)] \Lambda([c]) {\Bbb Q}^{\Phi}(t)(dc)
  \end{equation} 
  for the typical cell distribution ${\Bbb Q}^\Phi(t)$. Observe next that the use of (\ref{STDscaling})
  gives us the important scaling relation
  \begin{equation}\label{ScalRelTC}
   {\Bbb Q}^{\Phi}(t) = [t^{-1}] \odot {\Bbb Q}^{\Phi}(1),
  \end{equation}
  whence the desired equation (\ref{TypDet}) follows upon putting ${\Bbb Q}^\Phi:={\Bbb Q}^\Phi(1)$.\\
  To proceed, we let now $\tilde{\Bbb Q}^{\Phi}(t)$ stand for the image of
  the typical cell's law ${\Bbb Q}^{\Phi}(t)$ under the standard re-normalizing map ${\cal K}_0 \ni c 
  \mapsto \tilde{c} = \frac{1}{\Lambda([c])} c \in {\cal K}_{0,1}.$ 
  Hence
 \begin{equation}\label{QRenormConst}
  \tilde{\Bbb Q}^{\Phi}(t) = \tilde{\Bbb Q}^{\Phi} = {\rm const}.
 \end{equation}
 Consequently, (\ref{TypCellEqn}) leads us to
 $$ 0 = \int_{{\cal K}_{0,1}} \left(
  2 \Phi^{\oslash}_{0,1}(\cdot | c) - \delta_c - \tilde{\Bbb Q}^{\Phi}
  \right) d \tilde{\Bbb Q}^{\Phi}(dc) $$
 and thus 
 $$ \tilde{\Bbb Q}^{\Phi} = \int_{{\cal K}_{0,1}} \Phi^{\oslash}_{0,1}(\cdot|c) 
     \tilde{\Bbb Q}^{\Phi}(dc), $$
 so that eventually, but not unexpectedly, $\tilde{\Bbb Q}^{\Phi}$ turns out to coincide
 with the unique invariant measure $\varpi_{\Phi}$ for the re-normalized kernel 
 $\Phi^{\oslash}_{0,1}$, i.e.
 \begin{equation}\label{QRenorm}
  \tilde{\Bbb Q}^{\Phi} = \varpi_{\Phi}
 \end{equation}
 which is claimed as (\ref{TypCellRel1}) in (d).\\ Getting back to the original non-rescaled typical cell ${\Bbb Q}^{\Phi}$ constituting the focus
 of our present interest, from (\ref{QRenorm}) we obtain the integral representation
 \begin{equation}\label{RepreQ}
  {\Bbb Q}^{\Phi}(t) = \int_{{\Bbb R}_+} \left( s \odot \varpi_{\Phi} \right) \nu(ds|t),
 \end{equation}
 where $\nu(ds|t) = {\Bbb Q}^{\Phi}(t)(\left\{\Lambda([c]) \in ds\right\}) = 
 {\Bbb P}(\{\Lambda([\TypicalCell(Y^{\Phi}(t))]) \in ds\})$ is the law of the typical cell's split
 intensity/mean width under ${\Bbb Q}^{\Phi}(t).$ However, (\ref{EqDistr}) is just a reformulation of (\ref{RepreQ}).\\ The corresponding statement from part (f) for the zero cell of the tessellation follows now readily by the general theory, see e.g. \cite[Thm 10.4.1]{SW}, because of the finite mean volume of the typical cell, which was part of statement (d) of the main Theorem.\\
 Further, in part (e) the claimed preservation of ${\Bbb Q}^{\Phi}$ by the {\bf [Continuous Shrink Dynamics]} follows 
now by
\begin{itemize}
 \item
 using (\ref{ScalRelTC}) to write $[1+dt] \odot {\Bbb Q}^{\Phi} = {\Bbb Q}^{\Phi}(1-dt)$
and comparing with {\bf [CSD1]},
\item 
 noting thereupon that during the time interval $[1-dt,1]$ of the incremental MNW-construction,
 when constructing $Y^{\Phi}(1)$ out from $Y^{\Phi}(1-dt),$ precisely the contents of rule
 {\bf [CSD2]} is carried out for each cell of $Y^{\Phi}(1-dt).$
\end{itemize}
 Alternatively one can just use (\ref{TypCellEqn}) at $t=1.$ Finally, the required uniqueness of the stationary regime of the {\bf [Continuous Shrink Dynamics]} follows from 
(\ref{FstOrd}) together with the observation that under each of the conditions ${\Bbb E}\Lambda([\hat{K}_t^{\oslash;0}])=\alpha>0$ there is precisely one stationary regime by assumption {\bf [SKR2]}.

\paragraph{Identifying the translation invariant field as the outcome of the global construction.}
  To proceed with our argument, consider the sequence $(S^*_i,\tau^*_i)_{i \geq 0},\;
  S^*_{i+1} \supset S^*_i,$ of successive time-marked cells of $Y^{\Phi}(t)$
  containing the origin ${\bf 0}.$ Then it is clear that $(S^*_i,\tau^*_i)_{i \geq 0}$
  is a realization of time-marked spinal growth process backwards in time as
  determined by the extended spinal growth kernel $\hat{\Phi}^*$ given by (\ref{JointTMGrow})
 and with initial seed cell $S^*_0 = S_{\rm seed}$ given in law as the zero cell (\ref{ZCell}) of $Y^\Phi(t)$.
 Thus, $Y^{\Phi}(t)$ admits a spinal frame arising as an instance of {\bf [Construction: Spinal
 phase]} initiated by the zero cell (\ref{ZCell}). Moreover, once this spinal frame for $Y^{\Phi}(t)$ 
 is given, the required compatibility of $Y^{\Phi}(t)$ immediately implies that the remaining
 part of the tessellation has to coincide with the outcome of an instance of our {\bf [Construction:
 Local fill phase]} as required.

\paragraph{Uniqueness in the translation invariant regime.}
 This follows directly by the previous paragraph where our identification of $Y^{\Phi}(t)$ with
 the outcome of global construction was only reliant on the fact that $Y^{\Phi}(t)$ is $\Phi$-compatible.

\paragraph{Uniqueness in the general regime.}
 Write  ${\Bbb Z}^{\Phi}_1$ for the zero-cell of $Y^{\Phi}(t)$ and set $\zeta^{\Phi}_1$ to be
 its corresponding birth time in the course of the recursive split construction. Further, let recursively 
 ${\Bbb Z}^{\Phi}_{k+1}$ be the parent cell of ${\Bbb Z}^{\Phi}_k$ and, again, $\zeta_{k+1}^\Phi$
 its corresponding  birth time. Next, consider some arbitrary $\Phi$-compatible nested tessellation $Y'(t)$ and define $({\Bbb Z}'_k,\zeta'_k)$ for $Y'(t)$ in analogy to 
$({\Bbb Z}^{\Phi}_k,\zeta^{\Phi}_k)$ for $Y^{\Phi}(t).$ Clearly, both $({\Bbb Z}'_k,\zeta'_k)$
 and $({\Bbb Z}^{\Phi}_k,\zeta^{\Phi}_k)$ are Markov processes governed by the same Markovian
 dynamics driven by suitable growth kernels, recall (\ref{JointTMGrow}) in Subsection \ref{GlConstr}.
 In view of the assumed \textbf{[Split kernel regularity]} we conclude that
 \begin{equation}\nonumber
  \lim_{k \to \infty} {\rm d}_{{\rm TVar}}\left( ({\Bbb Z}'_k,\zeta'_k), 
  ({\Bbb Z}^{\Phi}_k,\zeta^{\Phi}_k) \right) = 0,
 \end{equation}
 with ${\rm d}_{{\rm TVar}}$ standing for the total variation distance of distributions, cf. \cite[Thm. 20.20]{KB}.
 Consequently, given $\epsilon > 0$ the processes $Y'(t)$ and $Y^{\Phi}(t)$ can be coupled
 so that with probability exceeding $1-\epsilon$ we have 
 ${\Bbb Z}^{\Phi}_k = {\Bbb Z}'_k$ and $\zeta^{\Phi}_k = \zeta'_k$ as soon as $k$ is large enough,
 say $k \geq k[\epsilon].$ We modify this coupling by erasing, on the above event of equality, 
 both $Y'(t)$ and $Y^{\Phi}(t)$ over the common region ${\Bbb Z}^{\Phi}_k = {\Bbb Z}'_k,$
 and then restoring both over this region using the same instance of the recursive split
 construction. The resulting new coupling has the property that, with probability not smaller than
 $1-\epsilon,$ we get the equality $Y'(t) \cap {\Bbb Z}'_k = Y^{\Phi}(t) \cap {\Bbb Z}'_k$
 for $k$ large enough. Recalling that almost surely ${\Bbb Z}'_k \uparrow {\Bbb R}^d$
 shows that ${\rm d}_{{\rm TVar}}(Y'(t),Y^{\Phi}(t)) \leq
 \epsilon.$ Taking into account that $\epsilon > 0$ was arbitrarily small we finally conclude that
 $Y'(t)$ and $Y^{\Phi}(t)$ coincide in law, as required.\hfill $\Box$

\begin{appendix}
\section{Argument for Proposition \ref{PHIINF}}\label{SecAp}
In order to keep the paper self-contained, we sketch here the proof of Proposition \ref{PHIINF} and refer to Lemma 1 of \cite{SLDP} for more details in a more general context.\\ \\ We start by defining for $y\in{\Bbb Y}_{[0,t]}(W)$ the split measure on $[W]\times W$ by $$\Upsilon_{\Phi_s^W}^y:=\sum_{c\in\Cells(y)}\int_{[c]}\delta_{(H,\pi_{\bf 0}(c\cap H))}\Phi_s^W(dH|c,y),$$ where $\pi_{\bf 0}$ stands for the metric projection of the origin $\bf 0$ onto the argument closed convex set and $\delta_{(\cdot)}$ denotes the unit mass Dirac measure. Write $\Facet_s(x,H)$ for the facet that arises by splitting by $H$ at time $s$ the unique cell $c\in\Cells(y)$ with $\pi_{\bf 0}(c\cap H)=x$ if such a facet exists and put $\Facet_s(x,H):=\emptyset$ otherwise. Then the generator ${\Bbb L}_s^W$ can be expressed as $$[{\Bbb L}_s^Wf](y)=\int_{[W]\times W}[f(y\cup\Facet_s(x,H))-f(y)]\Upsilon_{\Phi_s^W}^y(dH,dx).$$ Taking into account that $Y(s,V)$ arises as restriction of $Y(s,W)$ to $V$ for all compact and convex $V\subset W$, we get $$[{\Bbb L}_s^Vf](y)={\Bbb E}[[{\Bbb L}_s^Vf](Y(s,W))|Y(s,W)\cap V=y].$$ Let $\Im_s^W$ be the $\sigma$-field generated by $Y^\Phi(s)$, $W\subseteq{\Bbb R}^d$, $s\in[0,t]$. Defining the $\Im_s^W$-measurable random measure $\Upsilon_{\Phi_s^W}$ by $\Upsilon_{\Phi_s^W}:=\Upsilon_{\Phi_s^W}^{Y(s,W)}$, we conclude from our discussion above that for any measurable $B\subset[V]\times V$ we have $$\Upsilon_{\Phi_s^V}(B)={\Bbb E}[\Upsilon_{\Phi_s^W}(B)|{\Im}_s^V].$$ Noting that the measures $\Upsilon_{\Phi_s^W}$ are non-negative, the martingale convergence theorem ensures the existence of an ${\Im}_s^{{\Bbb R}^d}$-measurable random measure $\Upsilon_{\Phi_s}$ on ${\cal H}\times{\Bbb R}^d$ satisfying $$\Upsilon_{\Phi_s^W}(B)={\Bbb E}[\Upsilon_{\Phi_s}(B)|{\Im}_s^W]$$ for $B$ as above and any compact convex $W\subset{\Bbb R}^d$. Since ${\Bbb Y}_{[0,s]}$ is a Polish, there exists a measurable map $\Upsilon_{\Phi_s}$ from ${\Bbb Y}_{[0,s]}$ to the space of non-negative Borel measures on ${\cal H}\times{\Bbb R}^d$ with $\Upsilon_{\Phi_s}=\Upsilon_{\Phi_s}^{Y(s)}$ and we readily conclude that the generator of the (non-homogeneous) Markov process $Y(s)$, $s\in[0,t]$ is given by $$[{\Bbb L}_sf](y)=\sum_{c\in\Cells(y)}\int_{[y]}[f(y[c\oslash_s H])-f(y)]\Phi_s(dH|c,y),$$ where the last expression is defined for all bounded $f$, which are ${\Im}_t^W$-measurable for some bounded convex $W\subset{\Bbb R}^d$. The remaining relation (\ref{ZalMdzygen}) is an easy consequence of our above argument.\hfill $\Box$
\end{appendix}

\begin{flushleft}
\small
Tomasz Schreiber\hfill Christoph Th\"ale\\
Faculty of Mathematics and Computer Science\hfill Institute of Mathematics\\
Nicolaus Copernicus University\hfill University of Osnabr\"uck\\
Toru\'n, Poland\hfill Osnabr\"uck, Germany\\
\hfill {\tt christoph.thaele[at]uni-osnabrueck.de}
\end{flushleft}

\end{document}